\documentstyle[amssymb,12pt]{amsart}
\newtheorem{prop}{Proposition}[section]
\newtheorem{prop:def}{Proposition-Definition}[section]

\newtheorem{lemma}{Lemma}[section]
\newtheorem{thm}{Theorem}[section]
\newtheorem{cor}{Corollary}[section]

\theoremstyle{remark}
\newtheorem{remark}{Remark}

\begin{document} 

\newcommand{\nc}{\newcommand} \nc{\on}{\operatorname}

\nc{\¦}{{|}}

\nc{\pa}{\partial} 

\nc{\cA}{{\cal A}}\nc{\cB}{{\cal B}}\nc{\cC}{{\cal C}} \nc{\cE}{{\cal E}}
\nc{\cG}{{\cal G}}\nc{\cH}{{\cal H}}\nc{\cX}{{\cal X}}\nc{\cR}{{\cal R}}
\nc{\cL}{{\cal L}} \nc{\cK}{{\cal K}}\nc{\cO}{{\cal O}} 
\nc{\cF}{{\cal F}}\nc{\cM}{{\cal M}} \nc{\cW}{{\cal W}} 

\nc{\sh}{\on{sh}}\nc{\Id}{\on{Id}}\nc{\Diff}{\on{Diff}}
\nc{\ad}{\on{ad}}\nc{\Der}{\on{Der}}\nc{\End}{\on{End}}
\nc{\res}{\on{res}}\nc{\ddiv}{\on{div}} \nc{\FS}{\on{FS}}
\nc{\card}{\on{card}}\nc{\dimm}{\on{dim}}\nc{\gr}{\on{gr}}
\nc{\Jac}{\on{Jac}}\nc{\Ker}{\on{Ker}} \nc{\Den}{\on{Den}}
\nc{\Imm}{\on{Im}}\nc{\limm}{\on{lim}}\nc{\Ad}{\on{Ad}}
\nc{\ev}{\on{ev}} \nc{\Hol}{\on{Hol}}\nc{\Det}{\on{Det}}
\nc{\Cone}{\on{Cone}} \nc{\pseudo}{{\on{pseudo}}}
\nc{\class}{{\on{class}}}\nc{\rat}{{\on{rat}}}
\nc{\local}{{\on{local}}}\nc{\an}{{\on{an}}}
\nc{\Lift}{{\on{Lift}}}\nc{\Mer}{{\on{Mer}}}\nc{\mer}{{\on{mer}}}
\nc{\lift}{{\on{lift}}}\nc{\diff}{{\on{diff}}}\nc{\Aut}{{\on{Aut}}}
\nc{\DO}{{\on{DO}}}\nc{\Frac}{{\on{Frac}}}\nc{\cl}{{\on{class}}}
\nc{\Fil}{{\on{Fil}}}

\nc{\Bun}{\on{Bun}}\nc{\diag}{\on{diag}}\nc{\KZ}{{\on{KZ}}}
\nc{\CB}{{\on{CB}}}\nc{\out}{{\on{out}}}\nc{\Hom}{{\on{Hom}}}
\nc{\FO}{{\on{FO}}}

\nc{\al}{\alpha}\nc{\de}{\delta}\nc{\si}{\sigma}\nc{\ve}{\varepsilon}\nc{\z}{\zeta}
\nc{\vp}{\varphi} \nc{\la}{{\lambda}}\nc{\g}{\gamma}\nc{\eps}{\epsilon}
\nc{\PsiDO}{\Psi\on{DO}}

\nc{\AAA}{{\mathbb A}}\nc{\CC}{{\mathbb C}}\nc{\NN}{{\mathbb N}}
\nc{\PP}{{\mathbb P}}\nc{\RR}{{\mathbb R}}\nc{\VV}{{\mathbb V}}
\nc{\ZZ}{{\mathbb Z}} 

\nc{\bla}{{\mathbf \lambda}} \nc{\bv}{{\mathbf v}}
\nc{\bz}{{\mathbf z}}\nc{\bt}{{\mathbf t}} 
\nc{\bP}{{\mathbf P}} \nc{\kk}{{\mathbf k}}

\nc{\A}{{\mathfrak a}}\nc{\B}{{\mathfrak b}}\nc{\G}{{\mathfrak g}}
\nc{\HH}{{\mathfrak h}}\nc{\mm}{{\mathfrak m}}\nc{\N}{{\mathfrak n}}
\nc{\SG}{{\mathfrak S}}\nc{\La}{\Lambda}

\nc{\wt}{\widetilde}
\nc{\wh}{\widehat} \nc{\bn}{\begin{equation}}\nc{\en}{\end{equation}}
\nc{\SL}{{\mathfrak{sl}}}\nc{\ttt}{{\mathfrak{t}}}
\nc{\GL}{{\mathfrak{gl}}}

\title[Commuting families in skew fields and Beauville's fibration]
{Commuting families in skew fields and  quantization of Beauville's 
fibration}

\begin{abstract}
We construct commuting families in fraction fields of symmetric powers
of algebras. The classical limit of this construction gives 
Poisson commuting families associated with linear systems. 
In the case of a K3 surface $S$, they correspond to lagrangian fibrations
introduced by Beauville. When $S$ is the canonical cone of an algebraic
curve $C$, we construct commuting families of differential operators on 
symmetric powers of $C$, quantizing the Beauville systems. 
\end{abstract}

\author{B. Enriquez and V. Rubtsov}

\address{B.E.: IRMA, Universit\'e  Louis Pasteur, 7, rue Ren\'e Descartes, 
67084 Strasbourg, France}

\address{V.R.: D\'epartement de Math\'ematiques, Universit\'e d'Angers, 
49045 Angers, France}
\maketitle

\subsection*{Introduction}

In \cite{Beauville}, Beauville introduced lagrangian fibrations 
associated with a K3 surface $S$. These fibrations have the 
form $S^{[g]} \to \PP(H^0(S,\cL))$, where $S^{[g]}$ is the Hilbert scheme of 
$g$ points of $S$, equipped with a symplectic structure introduced in \cite{Mukai}, 
and $\cL$ is a line bundle on $S$.  Later, the authors of \cite{Donagi}
explained that these systems are natural deformations of 
the "separated" (in the sense of \cite{GNR}) versions of Hitchin's integrable systems, 
more precisely, of their description in terms of spectral curves
(already present in \cite{Hitchin}). Beauville's systems can be 
generalized to surfaces with a Poisson structure (see \cite{Bottacin}).
When $S$ is the canonical cone $\Cone(C)$ of an algebraic curve $C$, then 
this system coincides with the separated version of Hitchin's systems. 

A quantization of Hitchin's system was proposed in \cite{BD}. 
It seems interesting to construct quantizations of Beauville's 
systems. 

In this paper, we construct generalizations of the birational 
version of Beauville's construction (Theorem \ref{thm:Poisson})
and a quantum analogue of this construction (Theorem \ref{thm:comm}). 
We show that to obtain a quantization of Beauville's fibration, 
it would be sufficient to have quantizations of function fields of K3 surfaces. 
Such a quantization is not known explicitly, in general. 

However, in the case of the canonical cone of an algebraic curve $C$, an explicit 
quantization is known (see \cite{CMZ,EO}). In Section \ref{sect:final}, 
we construct the quantized Beauville systems explicitly in this case, and 
we show that in some cases these systems correspond to 
commuting families of rational differential operators on symmetric powers 
of $C$. We make these operators explicit in the case of a rational 
curve with marked points. 

Finally, in Appendix \ref{app}, we discuss the relation of 
Theorems \ref{thm:comm} and \ref{thm:Poisson} with the formal 
non-commutative geometry of \cite{Kap}.

\section{Commuting families in skew fields}

Let $A$ be an algebra with unit. For $f_1,\ldots,f_n$ elements of $A$, we set 
$$
[f_1,\ldots,f_n] = \sum_{\sigma\in \SG_n}
\eps(\sigma) f_{\sigma(1)} \otimes \cdots \otimes f_{\sigma(n)} ;   
$$ 
so $[f_1,\ldots,f_n]$ belongs to $A^{\otimes n}$. Recall that a skew field is 
an algebra with unit, where each nonzero element is invertible. 

\begin{thm} \label{thm:comm}
Assume that there exists a pair $(\phi,F)$ of a skew field $F$
and a ring injection $\phi : A^{\otimes n} \hookrightarrow F$. 
(We will identify elements of $A^{\otimes n}$ with their images by $\phi$.)
Assume that $f_0,\ldots,f_n$ are linearly independent elements of $A$.  
Set 
$$
\Delta_i = [f_0,\ldots,\check{f}_i, \ldots, f_n], 
$$ 
for $i = 0,\ldots,n$. Then $\Delta_0 \neq 0$. The elements 
$$
H_i = (\Delta_0)^{-1} \Delta_i
$$ 
form a commutative family of elements of $K$: we have 
$H_i H_j = H_j H_i$ for any $i,j$. 
\end{thm}

{\em Proof.} For $(i_1,\ldots,i_k)$ a family of distinct elements of 
$\{0,\ldots,n\}$ and $b\in A^{\otimes k}$, we define $b^{(i_1,\ldots,i_k)}$
as the image of $b$ by the injection $A^{\otimes k} \hookrightarrow A^{\otimes n}$, 
associated to $(i_1,\ldots,i_k)$. Since the family $(f_1,\ldots,f_n)$ is free, 
$\Delta_0 \neq 0$. Moreover, any subfamily of $(f_0,\ldots,f_n)$ is also 
free, so we can work by induction on $n$. 

We will prove the following identities
\begin{equation} \label{main:id}
\Delta_i (\Delta_0)^{-1} \Delta_j = \Delta_j (\Delta_0)^{-1} \Delta_i,  
\end{equation}
where $\Delta_i = [f_0,\ldots,\check f_i,\ldots,f_n]$; 
\begin{equation} \label{2a}
\sum_{i=1}^n (-1)^{i} [f_1,\ldots,\check f_i,\ldots,f_n]^{(1,\ldots,n-1)}
[f_1,\ldots,f_n]^{-1} (f_i)^{(n)} = (-1)^n , 
\end{equation}
and 
\begin{equation} \label{2b}
\sum_{i=1}^n (-1)^{i} [f_1,\ldots,\check f_i,\ldots,f_n]^{(1,\ldots,n-1)}
[f_1,\ldots,f_n]^{-1} (f_i)^{(a)} = 0  
\end{equation}
for $a = 1,\ldots,n-2$. Denote these identities by 
$(\ref{main:id}_n)$, $(\ref{2a}_n)$ and $(\ref{2b}_n)$. We will 
prove the implications 
$$
(\ref{2a}_2,\ref{2b}_2) \Rightarrow 
(\ref{main:id}_2) \Rightarrow 
(\ref{2a}_3,\ref{2b}_3) \Rightarrow 
(\ref{main:id}_3) \Rightarrow 
\cdots \Rightarrow 
(\ref{main:id}_{n-1}) \Rightarrow 
(\ref{2a}_n,\ref{2b}_n) \Rightarrow 
(\ref{main:id}_n) \Rightarrow 
\cdots
$$

Let us prove $(\ref{2a}_2)$.  This is the identity 
\begin{equation} \label{2a:part}
f^{(1)} [f,g]^{-1} g^{(2)} - g^{(1)} [f,g]^{-1} f^{(2)} = 1 
\end{equation}
(recall that $[f,g] = f^{(1)}g^{(2)} - g^{(1)} f^{(2)}$).  
We have 
$$
f^{(1)} [f,g]^{-1} g^{(2)} = \big( 1 - (gf^{-1})^{(1)}
(g^{-1}f)^{(2)} \big)^{-1} = (1-x)^{-1} , 
$$
where $x = gf^{-1} \otimes g^{-1}f$. On the other hand, we have 
$$
g^{(1)} [f,g]^{-1} f^{(2)} = \big( (fg^{-1})^{(1)} (f^{-1}g)^{(2)} - 1\big)^{-1} 
= (x^{-1} -1)^{-1} .  
$$
Then $(1-x)^{-1} + (x^{-1}-1)^{-1} = (1-x)^{-1}(1-x) = 1$, which
proves (\ref{2a:part}). 

Let us prove $(\ref{2b}_2)$.  This is the identity 
\begin{equation} \label{2b:part}
f^{(1)} [f,g]^{-1} g^{(1)} - g^{(1)} [f,g]^{-1} f^{(1)} = 0. 
\end{equation}
We have 
$$
f^{(1)} [f,g]^{-1} g^{(1)} = \big( (g^{-1})^{(1)} g^{(2)}
- (f^{-1})^{(1)} f^{(2)}\big)^{-1}
$$
and 
$$
g^{(1)} [f,g]^{-1} f^{(1)} = \big( (g^{-1})^{(1)} g^{(2)}
- (f^{-1})^{(1)} f^{(2)}\big)^{-1}, 
$$
which proves (\ref{2b:part}). 

Let us prove $(\ref{main:id}_2)$. This is the identity 
\begin{align} \label{main:id:part}
& \big( f^{(1)} h^{(2)} - f^{(2)} h^{(1)} \big) 
\big( f^{(1)} g^{(2)} - f^{(2)} g^{(1)} \big)^{-1} 
\big( g^{(1)} h^{(2)} - g^{(2)} h^{(1)} \big) 
\\ & \nonumber 
= \big( g^{(1)} h^{(2)} - g^{(2)} h^{(1)} \big) 
\big( f^{(1)} g^{(2)} - f^{(2)} g^{(1)} \big)^{-1} 
\big( f^{(1)} h^{(2)} - f^{(2)} h^{(1)} \big) . 
\end{align}
Write the difference of both sides as 
\begin{equation} \label{star}
h^{(1)} A h^{(1)} +  h^{(1)} B h^{(2)} +  
h^{(2)} C h^{(1)} +  h^{(2)} D h^{(2)},  
\end{equation}
where 
$$
A = f^{(2)} \big( f^{(1)} g^{(2)} - f^{(2)} g^{(1)} \big)^{-1} g^{(2)}
- g^{(2)} \big( f^{(1)} g^{(2)} - f^{(2)} g^{(1)} \big)^{-1} f^{(2)}, 
$$
$$
B = - f^{(2)} \big( f^{(1)} g^{(2)} - f^{(2)} g^{(1)} \big)^{-1} g^{(1)}
+ g^{(2)} \big( f^{(1)} g^{(2)} - f^{(2)} g^{(1)} \big)^{-1} f^{(1)}, 
$$
$$
C = -f^{(1)} \big( f^{(1)} g^{(2)} - f^{(2)} g^{(1)} \big)^{-1} g^{(2)}
+ g^{(1)} \big( f^{(1)} g^{(2)} - f^{(2)} g^{(1)} \big)^{-1} f^{(2)}, 
$$
$$
D = f^{(1)} \big( f^{(1)} g^{(2)} - f^{(2)} g^{(1)} \big)^{-1} g^{(1)}
- g^{(1)} \big( f^{(1)} g^{(2)} - f^{(2)} g^{(1)} \big)^{-1} f^{(1)}.  
$$
As we have seen, $A = D= 0$, $B=1$, $C = -1$, so (\ref{star}) is equal to zero.  
This proves (\ref{main:id:part}). 

Let us now assume that the identities $(\ref{2a}_k)$, $(\ref{2b}_k)$ and 
$(\ref{main:id}_k)$ are proved for $k<n$. Let us prove $(\ref{2a}_n)$. 
Set
$$
T_i = (-1)^i [f_1,\ldots,\check f_i,\ldots,f_n]^{(1,\ldots,n-1)}
[f_1,\ldots,f_n]^{-1} (f_i)^{(n)}. 
$$ 
We want to prove 
\begin{equation} \label{2a:bis}
\sum_{i=1}^n T_i = (-1)^{n}. 
\end{equation}
Recall that 
$$
[f_1,\ldots,f_n] = \sum_{j=1}^n (-1)^{j+n+1} (f_j)^{(n)}
[f_1,\ldots,\check f_j,\ldots,f_n]^{(1,\ldots,n-1)}. 
$$
So 
\begin{align*}
& T_i = 
\\ & 
= \Big( \sum_{j=1}^n (-1)^{j+n+1} (f_i^{-1}f_j)^{(n)}
\big( [f_1,\ldots,\check f_j,\ldots,f_n] 
[f_1,\ldots,\check f_i,\ldots,f_n]^{-1}
\big)^{(1,\ldots,n-1)} \Big)^{-1} 
\\ & 
= \Big( \sum_{j=1}^n (-1)^{j+n+1} (f_i^{-1}f_1)^{(n)} (f_1^{-1}f_j)^{(n)}
\\ & 
\big( [f_1,\ldots,\check f_j,\ldots,f_n]
[\check f_1,\ldots,f_n]^{-1} 
[\check f_1,\ldots,f_n] 
 [f_1,\ldots,\check f_i,\ldots,f_n]^{-1}
\big)^{(1,\ldots,n-1)} \Big)^{-1} ; 
\end{align*}
($\ref{main:id}_{n-1}$) then allows to permute 
$$
[f_1,\ldots,\check f_j,\ldots,f_n][\check f_1,\ldots,f_n]^{-1}
\; \on{and} \; 
[\check f_1,\ldots,f_n][f_1,\ldots,\check f_i,\ldots,f_n]^{-1}.
$$ 
Then we get 
$$
T_i = (-1)^{i+1} T_1 (f_1^{-1}f_i)^{(n)} 
\big( 
[f_1,\ldots,\check f_i,\ldots,f_n][\check f_1,\ldots,f_n]^{-1}
\big)^{(1,\ldots,n-1)} . 
$$
Summing up, we get 
\begin{align*}
& \sum_{i=1}T_i 
= T_1 \sum_{i=1}^n (-1)^{i+1} \big((f_1)^{-1} f_i\big)^{(n)}
\big( [f_1,\ldots,\check f_i,\ldots,f_n]
[\check f_1,\ldots,f_n]^{-1} \big)^{(1,\ldots,n-1)}
\\ & = 
T_1 \big( (f_1)^{-1} \big)^{(n)} 
\big( \sum_{i=1}^n (-1)^{i+1} 
[f_1,\ldots,\check f_i,\ldots,f_n]^{(1,\ldots,n-1)} (f_i)^{(n)} \big) 
\big( [\check f_1,\ldots,f_n]^{-1} \big)^{(1,\ldots,n-1)}
\\ & 
= 
T_1 \big( (f_1)^{-1} \big)^{(n)}
(-1)^{n+1} [f_1,\ldots,f_n]
\big( [\check f_1,\ldots,f_n]^{-1} \big)^{(1,\ldots,n-1)}
\\ & = (-1)^{n}. 
\end{align*}
This proves ($\ref{main:id}_n$). 

Assume now that we have proved the identities $(\ref{2a}_k)$, $(\ref{2b}_k)$ 
and $(\ref{main:id}_k)$ for $k<n$, and let us prove $(\ref{2a}_n)$. 
This identity is equivalent to 
$$
\sum_{i=1}^n (-1)^{i} [f_1,\ldots,\check f_i,\ldots,f_n]^{(1,\ldots,n-1)}
[f_1,\ldots,f_n]^{-1} (f_i)^{(1)} = 0.  
$$
Let us set 
$$
U_i = (-1)^{i} [f_1,\ldots,\check f_i,\ldots,f_n]^{(1,\ldots,n-1)}
[f_1,\ldots,f_n]^{-1} (f_i)^{(1)}.  
$$
Then as in the proof of ($\ref{2a}_n$) (in particular, using ($\ref{main:id}_{n-1}$)),
one shows that 
$$
U_i = U_1 \cdot (-1)^{i+1} [f_1,\ldots,\check f_i,\ldots,f_n]^{(1,\ldots,n-1)}
\big( [f_2,\ldots,f_n]^{-1}\big)^{(1,\ldots,n-1)} (f_i)^{(1)}. 
$$ 
So ($\ref{2a}_n$) will be a consequence of 
\begin{equation} \label{3}
(f_1)^{(1)} + \sum_{i=2}^n (-1)^{i+1} [f_1,\ldots,\check f_i,\ldots,f_n]
[f_2,\ldots,f_n]^{-1} (f_i)^{(1)} = 0.  
\end{equation}
Let us prove (\ref{3}). By ($\ref{2a}_{n-1}$) and ($\ref{2b}_{n-1}$), we have 
\begin{equation} \label{USAF}
1 + \sum_{i=2}^n (-1)^{i+1} [f_2,\ldots,\check f_i ,\ldots,f_n]^{(2,\ldots,n)}
[f_2,\ldots,f_n]^{-1} (f_i)^{(1)} = 0, 
\end{equation}
and for $j = 2,\ldots,n$, we have 
\begin{equation} \label{B52}
\sum_{i=2}^n (-1)^{i+1} 
[f_2,\ldots,\check f_i,\ldots,f_n]^{(1,\ldots,\check j,\ldots,n)}
[f_2,\ldots,f_n]^{-1} (f_i)^{(1)} = 0, 
\end{equation}
so multiplying (\ref{USAF}) by $(f_1)^{(1)}$ and (\ref{B52}) by 
$(-1)^{j+1}(f_1)^{(j)}$ from the left, and adding up the resulting 
identities, we obtain (\ref{3}). 

Let us now assume that the identities ($\ref{main:id}_{n'}$), 
($\ref{2a}_{n'}$) and ($\ref{2b}_{n'}$) are proved for $n'<n$, as well as 
($\ref{2a}_{n}$) and ($\ref{2b}_{n}$), and let us prove ($\ref{main:id}_{n}$). 
Since we have 
$$
[f_1,\ldots,f_n] = \sum_{k=1}^n (-1)^{k+1} (f_1)^{(k)}
[f_2,\ldots,f_n]^{(1,\ldots,\check k,\ldots,n)}, 
$$
$\Delta_i (\Delta_0)^{-1} \Delta_j$ has the expansion
\begin{align*}
& \Delta_i (\Delta_0)^{-1} \Delta_j = 
\\ & 
\sum_{k,\ell = 1}^n (-1)^{k+\ell}
(f_0)^{(k)} [f_0,\ldots,\check f_i,\ldots,f_n]^{(1,\ldots,\check k,\ldots,n)} 
(\Delta_0)^{-1} [f_0,\ldots,\check f_j,\ldots,f_n]^{(1,\ldots,\check \ell,\ldots,n)} 
(f_0)^{(\ell)}  
\end{align*}
To prove ($\ref{main:id}_n$), we will prove 
\begin{align} \label{cordovero}
& [f_0,\ldots,\check f_i,\ldots,f_n]^{(1,\ldots,\check k,\ldots,n)} 
(\Delta_0)^{-1} [f_0,\ldots,\check f_j,\ldots,f_n]^{(1,\ldots,\check k,\ldots,n)} 
\\ & \nonumber 
= [f_0,\ldots,\check f_j,\ldots,f_n]^{(1,\ldots,\check k,\ldots,n)} 
(\Delta_0)^{-1} [f_0,\ldots,\check f_i,\ldots,f_n]^{(1,\ldots,\check k,\ldots,n)} 
\end{align} 
for any $k = 1,\ldots,n$, 
and 
\begin{align} \label{luria}
& [f_0,\ldots,\check f_i,\ldots,f_n]^{(1,\ldots,\check k,\ldots,n)} 
(\Delta_0)^{-1} [f_0,\ldots,\check f_j,\ldots,f_n]^{(1,\ldots,\check \ell,\ldots,n)} 
\\ & \nonumber 
- [f_0,\ldots,\check f_j,\ldots,f_n]^{(1,\ldots,\check \ell,\ldots,n)} 
(\Delta_0)^{-1} [f_0,\ldots,\check f_i,\ldots,f_n]^{(1,\ldots,\check k,\ldots,n)} 
\\ & \nonumber 
= (-1)^{i+j} [f_0,\ldots,\check f_i,\ldots,
\check f_j,\ldots,f_n]^{(1,\ldots,\check k,\ldots,\check \ell,\ldots,n)} 
\end{align} 
when $1\leq k < \ell \leq n$. We denote these identities by ($\ref{cordovero}_k$) 
and ($\ref{luria}_{k,\ell}$).  

Let us show why the collection of identities 
($\ref{cordovero}_k$) and ($\ref{luria}_{k,\ell}$) imply 
($\ref{main:id}_n$). We have 
\begin{align*}
& \Delta_i (\Delta_0)^{-1} \Delta_j - \Delta_j (\Delta_0)^{-1} \Delta_i 
\\ & =   
\sum_{k = 1}^n (f_0)^{(k)} (\on{l.h.s.\  --\  r.h.s.\ of\ }(\ref{cordovero}_{k})) 
(f_0)^{(k)}
\\ & 
+\sum_{1\leq k< \ell \leq n}
(f_0)^{(k)} (\on{l.h.s.\ of\ }(\ref{luria}_{k,\ell})) 
(f_0)^{(\ell)}
- (f_0)^{(\ell)} (\on{l.h.s.\ of\ }(\ref{luria}_{k,\ell})) 
(f_0)^{(k)}. 
\end{align*}
Then ($\ref{cordovero}_k$) implies that the first summand is 
zero, and ($\ref{luria}_{k,\ell}$) implies that the second summand
is zero, so $\Delta_i (\Delta_0)^{-1} \Delta_j 
- \Delta_j (\Delta_0)^{-1} \Delta_i = 0$, that is ($\ref{main:id}_n$). 

Let us now show why the identities ($\ref{main:id}_{n'}$), ($\ref{2a}_{n'}$), 
($\ref{2b}_{n'}$), $n'<n$, together with ($\ref{2a}_{n}$), 
($\ref{2b}_{n}$), imply the identities ($\ref{cordovero}_k$) and 
($\ref{luria}_{k,\ell}$).

Let us first prove ($\ref{cordovero}_k$). 
($\ref{main:id}_{n-1}$)  implies the identities 
\begin{align} \label{bernheim}
& [f_1,\ldots,\check f_j,\ldots,f_n]^{-1}
[f_1,\ldots,\check f_\ell,\ldots,f_n]
[f_1,\ldots,\check f_i,\ldots,f_n]^{-1}
\\& \nonumber 
= [f_1,\ldots,\check f_i,\ldots,f_n]^{-1}
[f_1,\ldots,\check f_\ell,\ldots,f_n]
[f_1,\ldots,\check f_j,\ldots,f_n]^{-1}. 
\end{align}
Insert each identity (\ref{bernheim}) in the factors 
$(1,\ldots,\check k,\ldots,n)$, multiply it by 
$(-1)^\ell (f_\ell)^{(k)}$, and sum up all resulting identities. 
Using the expansion of $\Delta_0$ in the form 
$$
\Delta_0 = (-1)^k \sum_{\ell = 1}^{n} (-1)^\ell (f_\ell)^{(k)}
[f_1,\ldots,\check f_{\ell},\cdots f_n]^{(1,\ldots,\check k,\ldots,n)} ,  
$$
this identity can then be written as follows
\begin{align*} 
& \big( [f_0,\ldots,\check f_i,\ldots,f_n]^{-1} \big)^{(1,\ldots,\check k,\ldots,n)} 
\Delta_0 \big( [f_0,\ldots,\check f_j,\ldots,f_n] \big)^{(1,\ldots,\check k,\ldots,n)} 
\\ & 
= \big( [f_0,\ldots,\check f_j,\ldots,f_n] \big)^{(1,\ldots,\check k,\ldots,n)} 
\Delta_0 \big( [f_0,\ldots,\check f_i,\ldots,f_n] \big)^{(1,\ldots,\check k,\ldots,n)}.  
\end{align*} 
Taking the inverse of this identity, we obtain the identity ($\ref{cordovero}_k$).

Let us now prove identity ($\ref{luria}_{k,\ell}$). Using the symmetries of the 
brackets $[f_1,\ldots,f_n]$, this identity can be written 
as follows 
\begin{align} \label{2prime}
& [f_1,\ldots,f_{n-2},f_{n-1}]^{(1,\ldots,n-2,n-1)}
[f_1,\ldots,f_{n-2},f_{n-1},f_n]^{-1} 
[f_1,\ldots,f_{n-2},f_n]^{(1,\ldots,n-2,n)}
\\& \nonumber 
- [f_1,\ldots,f_{n-2},f_n]^{(1,\ldots,n-2,n-1)}
[f_1,\ldots,f_{n-2},f_{n-1},f_n]^{-1} 
[f_1,\ldots,f_{n-2},f_{n-1}]^{(1,\ldots,n-2,n)}
\\ & \nonumber
= [f_1,\ldots,f_{n-2}]^{(1,\ldots,n-2)}. 
\end{align}
Let us now prove (\ref{2prime}). Recall that we assume that 
the identities ($\ref{2a}_n$) and ($\ref{2b}_n$) are proved. 
Denote the identity ($\ref{2b}_n$) corresponding to index $a$, 
by ($\ref{2b}_{n,a}$). Multiply  ($\ref{2a}_n$) by 
$(-1)^n [f_1,\ldots,f_{n-2}]^{(1,\ldots,n-2)}$ from the right, 
multiply ($\ref{2b}_{n,a}$) by  $(-1)^a [f_1,\ldots,f_{n-2}]^{(1,\ldots,\check a, 
\ldots, n-2,n)}$ from the right, and sum up all resulting identities. 
Then the identities 
$$
\sum_{a=1}^{n-2} (-1)^a (f_i)^{(a)} [f_1,\ldots,\check f_i,\ldots,f_{n-2}
]^{(1,\ldots,\check a,\ldots,n-2)} = (-1)^i [f_1,\ldots,f_{n-2}]
$$
and $[f_1,\ldots,f_{n-2}] = 0$ when $f_i = f_j$ and $i\neq j$, yield
(\ref{2prime}). This proves that ($\ref{2a}_n$) and ($\ref{2b}_n$) 
imply ($\ref{main:id}_n$). 
\hfill \qed \medskip 

The proof of Theorem \ref{thm:comm} can be generalized to give the 
following result. 

\begin{cor}\label{comm:fij}
Let $\cA$ be an algebra, 
$(f_{i,j})_{0\leq i \leq n, 1\leq j \leq n}$ be elements of $\cA$ such that 
$$
f_{i,j}f_{k,\ell} = f_{k,\ell} f_{i,j}
$$
for any $i,j,k,\ell$ such that $j\neq\ell$. For any $I \subset 
\{0,\ldots,n\},J\subset \{1,\ldots,n\}$ of the same cardinality, we set 
$\Delta_{I,J} = \sum_{\sigma\in \on{Bij}(I,J)} \eps(\sigma) f_{i,\sigma(i)}$. 
Assume that the $\Delta_{I,J}$ are all invertible. Set 
$\Delta_i = \Delta_{\{1,\ldots,n\},\{0,\ldots,\check i, \ldots,n\}}$. 
Then the 
$$
H_i = (\Delta_0)^{-1} \Delta_i
$$
all commute together. 
\end{cor}

(Theorem \ref{thm:comm} can be recovered 
in when $f_{i,j} = (f_i)^{(j)}$; the assumption on the $\Delta_{I,J}$
is a replacement of the assumption of freeness of the family 
$(f_0,\ldots,f_n)$.) 

\begin{remark}
Assume that $A$ has an involution. Then if $f_0,\ldots,f_n$
are self-adjoint, so is each $\Delta_i$. If $\Delta_0$ has a self-adjoint 
square root $(\Delta_0)^{1/2}$, then the family 
$\wt H_i = (\Delta_0)^{-1/2} \Delta_i (\Delta_0)^{-1/2}$ is a commuting 
family of self-adjoint operators.  
\end{remark}


\section{Poisson commuting families on symmetric powers}

We will fix a base field $\kk$ of characteristic $\neq 2$. 

\subsection{Poisson commuting families}

\begin{lemma} \label{prol:poisson}
If $B$ is an integral Poisson algebra, then there is a 
unique Poisson structure on $\Frac(B)$ extending the Poisson 
structure of $B$. 
\end{lemma}

{\em Proof.} This structure is uniquely defined by the relations 
$\{1/f,g\} = - \{f,g\} / f^2$, $\{1/f,1/g\} =\{f,g\} / (f^2 g^2)$. 
\hfill \qed \medskip 

Theorem \ref{thm:comm} has a Poisson counterpart. 

\begin{thm} \label{thm:Poisson}
Let $A$ be a Poisson algebra. Assume that $A^{\otimes n}$
is integral, and let $\pi : A^{\otimes n} \hookrightarrow 
\Frac(A^{\otimes n})$ be its injection in its fraction field. 
For each free family $(f_0,\ldots,f_n)$ of elements of $A$, 
we set 
$$
\Delta_i^\cl = [f_0,\ldots,\check f_i, \ldots, f_n]. 
$$ 
Then $\Delta_0^\cl$ is nonzero. Set $H_i^\cl = 
\Delta_i^\cl / \Delta_0^\cl$. Then the family 
$(H_i^\cl)_{i = 1,\ldots,n}$ is Poisson-commutative: 
$$
\{H_i^\cl,H_j^\cl\} = 0
$$
for any pair $(i,j)$.  
\end{thm}

{\em Proof}. We will give two proofs. 

\medskip \noindent {\em First proof.}
Theorem \ref{thm:comm} may be extended to the case where $\Frac(A^{\otimes n})
\hookrightarrow K$ is replaced by an injection $\Frac(A^{\otimes n})
\hookrightarrow R$, where $R$ is an algebra where the 
$[f_{i_1},\ldots,f_{i_k}]^{(\al_1,\ldots,\al_k)}$
($1\leq i_1 < \cdots < i_k \leq n$, $\al_1,\ldots,\al_k$
all distinct)  are all invertible. Then we apply this generalization of 
Theorem \ref{thm:comm} to the following algebras: the algebra $A$
of Theorem \ref{thm:comm} is the $\kk[\eps] / (\eps^2)$-algebra 
$A[\eps]/(\eps^2)$, equipped with the product $f \cdot g = 
f g + \eps \{f,g\}$; the algebra $R$ is the $\kk[\eps]/ (\eps^2)$-algebra 
$K[\eps] / (\eps^2)$, with a product defined in the same way. 

An element $f$ of $K[\eps]/(\eps^2)$ is invertible iff its reduction 
$f_0$ modulo $\eps$ is nonzero; on the other hand, if $f,g\in 
K[\eps]/(\eps^2)$, $fg = gf$ is equivalent to $\{f_0,g_0\} = 0$. 
So the elements $[f_{i_1},\ldots,f_{i_k}]^{(\al_{1}\ldots\al_{k})}$ 
are all invertible in $K[\eps]/(\eps^2)$; applying the generalization of 
Theorem \ref{thm:comm} to $(f_0,\ldots,f_n)$, we get $H_i H_j = H_j H_i$, 
so $\{H_i^\cl,H_j^\cl\} = 0$. 

\medskip \noindent {\em Second proof.}
We have to prove 
\begin{equation} \label{var:jac}
H_i^\cl \{H_j^\cl,H_k^\cl\} + \on{cyclic\ permutation\ in\ }(i,j,k) = 0. 
\end{equation}
We have 
$$
H_i^\cl = \sum_{p=1}^n \sum_{\al = 0}^n (-1)^{p+\al} (f_\al)^{(p)}
(M_{\al,i})^{(1\ldots\check p \ldots n)} ,  
$$
where
$$
M_{\al,i} = (-1)^{1_{\al<i}} [f_0\ldots \check f_\al \ldots \check f_i
\ldots f_n]
$$
if $\al \neq i$ (we set $1_{\al<i} = 1$ if $\al<i$ and $0$ otherwise)
and $M_{\al,i} = 0$ if $\al = i$. 
Now we have 
$$
\{H_i^\cl,H_j^\cl\} = \sum_{p=1}^n \sum_{\al,\beta = 0}^n 
(-1)^{\al+\beta} (\{f_\al,f_\beta\})^{(p)}
(M_{\al,i}M_{\beta,j} - M_{\al,j}M_{\beta,i})^{(1\ldots \check p \ldots n)}, 
$$
so identity (\ref{var:jac}) is a consequence of 
\begin{equation} \label{aux}
\forall (i,j,k,\al,\beta,\gamma), 
\; 
\sum_{\sigma\in \on{Perm}(i,j,k)}
\eps(\sigma) M_{\al,\sigma(i)} M_{\beta,\sigma(j)}
M_{\gamma,\sigma(k)} = 0. 
\end{equation} 
When $\card\{\al,\ldots,k\} = 3$, this identity follows from the antisymmetry
relation $M_{i,j} + M_{j,i} = 0$. When $\card\{\al,\ldots,k\} = 4$ (resp., 
$5,6$), it follows from the following Grassmann identities (to get (\ref{aux}), 
one should set $V = (A^{\otimes n})^{\oplus n}$ and $\La$ some partial 
determinant). Let $V$ be a vector space. Then 

-- if $\Lambda \in \wedge^2(V)$, and $a,b,c,d\in V$, then 
$$
\La(a,b) \La(c,d) - \La(a,c)\La(b,d) + \La(a,d) \La(b,c) = 0 ; 
$$

-- if $\La\in \wedge^3(V)$ and $a,b,c,b',c' \in V$, then 
\begin{align*}
& \La(b,c,c')\La(a,c,b')\La(b,b',c') 
+ \La(b,c,b')\La(c,b',c')\La(a,b,c')
\\ & - \La(b,c,b')\La(a,c,c')\La(b,b',c') 
- \La(b,c,c')\La(c,b',c')\La(a,b,b')
= 0; 
\end{align*} 

-- if $\La\in \wedge^4(V)$ and $a,b,c,a',b',c'\in V$, then 
\begin{align} \label{toprove}
& \nonumber \La(b,c,b',c')\La(a,c,a',c')\La(a,b,a',b') 
+ \La(b,c,a',c')\La(a,c,a',b')\La(a,b,b',c') 
\\ & \nonumber 
+ \La(b,c,a',b')\La(a,c,b',c')\La(a,b,a',c') 
\\ & \nonumber 
- \La(b,c,b',c')\La(a,c,a',b')\La(a,b,a',c') 
- \La(b,c,a',b')\La(a,c,a',c')\La(a,b,b',c') 
\\ & - \La(b,c,a',c')\La(a,c,b',c')\La(a,b,a',b') = 0.
\end{align}

Let us show these identities. Let $(\eps_1,\eps_2,\ldots)$ be a basis of 
$V$. In the first identity, we can assume by linearity that 
$\La = \eps_1 \wedge \eps_2$, then the identity is easy to check. 
In the second identity, we assume $\La = \eps_1 \wedge \eps_2 \wedge \eps_3$, 
then the identity takes place in a 3-dimensional vector space; 
then if rank$(b,c,b',c')<3$, the terms not involving $a$ are all 
zero, so the identity is satisfied; if rank$(b,c,b',c') = 3$, 
the identity is linear in $a$, so to check it, it is enough to 
check the identities where $a$ is replaced by $b,c,b',c'$, which 
are immediate. 

Let us prove identity (\ref{toprove}). We may assume that $\Lambda = 
\eps_1 \wedge \eps_2 \wedge \eps_3 \wedge \eps_4$, and then that $V$
is $4$-dimensional. Let $Q(a,b,c,a',b',c')$ be the l.h.s. of (\ref{toprove}). 
Then $Q$ is a polynomial on $V^6$. We may replace $\kk$ by its algebraic closure; 
then it is sufficient to prove the vanishing of $Q$ the open subset $U^2$ of all 
$(a,b,c,a',b',c')$, such that both families $(a,b,c)$ and $(a',b',c')$
are free. Moreover, we have $V^6 = (\kk^3 \oplus \kk^3) \otimes V$; 
there is an action of $SL_3(\kk) \otimes SL_3(\kk)$ on $V^6$, namely 
$(\rho\oplus \rho) \otimes \on{id}_V$, where $\rho$ is the vector representation of 
$SL_3(\kk)$ on $\kk^3$. Then one checks that if $(X,Y)\in \SL_3(\kk) \oplus 
\SL_3(\kk)$, we have 
$$
{d\over{d\eps}} Q((1+\eps X)(a,b,c),(1+\eps Y)(a',b',c'))_{\¦\eps = 0} = 0, 
$$
so $Q$ is constant along the orbits of $SL_3(\kk) \times SL_3(\kk)$ on $U^2$. 
These orbits are uniquely determined by the vector spaces spanned by 
$(a,b,c)$ and $(a',b',c')$, and by the volumes of $(a,b,c)$ and $(a',b',c')$
in these spaces. Since $\wedge^4(V)$ is $1$-dimensional, we may assume that 
$\Lambda = \eps'_1 \wedge \eps'_2 \wedge \eps'_3 \wedge \eps'_4$, where 
$(\eps'_1,\ldots,\eps'_4)$ is a basis of $V$ adapted to the subspaces 
$\on{Vect}(a,b,c)$ and $\on{Vect}(a',b',c')$. There are two possibilities: 

-- $\on{Vect}(a,b,c) = \on{Vect}(a',b',c')$. Then (\ref{toprove}) is trivially 
satisfied. 

-- dim$(\on{Vect}(a,b,c) \cap \on{Vect}(a',b',c')) = 2$. Then we may assume
$a = a' = \eps'_1$, $b = b' = \eps'_2$, $c = \la \eps'_3$, 
$c' = \mu \eps'_4$. Then  each triple product of (\ref{toprove}) vanishes, 
so (\ref{toprove}) also holds in this case. 
\hfill \qed \medskip

\subsection{The Beauville fibration associated to a K3 surface}
 
In this section, we explain, at the birational level, 
Beauville's construction of a lagrangian fibration associated to 
complex K3 surfaces (Proposition 3 of \cite{Beauville}). We show that 
this result can be rederived from Theorem \ref{thm:Poisson}.  

Let $S$ be a complex K3 surface, equipped with a very ample 
line bundle $\cL$. Set $g = h^0(S,\cL)-1$. Let $\varphi : S \hookrightarrow 
\PP(H^0(S,\cL)^*) = \CC P^g$ be the corresponding embedding. Let 
$\Cone(S)$ be the cone of this embedding, so we get an embedding 
$\Cone(S) \hookrightarrow \CC^{g+1}$ compatible with $\varphi$.

Let $A_S = \bigoplus_{i\geq 0} H^0(S,\cL^{\otimes i})$ be the graded 
algebra of $\Cone(S)$. Then $\varphi$ induces a morphism of graded algebras
$$
\varphi^* : \CC[X_0,\ldots,X_g] \to A_S. 
$$
Here $(X_0:\ldots:X_g)$ are projective coordinates on $\CC P^g$. 

If $B$ is any algebra, we set $B^{(g)} = (B^{\otimes g})^{\SG_g}$; 
then $\varphi^*$ induces an algebra morphism 
$$
(\varphi^*)^{(g)} : \CC[X_0,\ldots,X_g]^{(g)} \to (A_S)^{(g)}. 
$$

The geometric version of Beauville's construction (see \cite{Beauville}) 
is as follows. To a generic $g$-uple $(P_1,\ldots,P_g)$ of points of 
$\CC P^g$, we associate the unique hyperplane containing $(P_1,\ldots,P_g)$. 
The projective coordinates of this hyperplane are the minors of the matrix 
$(X_i^{(j)})_{0\leq i\leq n, 1\leq j\leq n}$, where for each 
$j$, $(X_i^{(j)})_{0\leq i\leq n}$ are projective coordinates of $P_j$. 
The affine coordinates of this hyperplane are therefore 
$((-1)^i h_i)_{i = 1,\ldots,g}$, where 
$$
h_i = \Delta_i(P_1,\ldots,P_g) / \Delta_0(P_1,\ldots,P_g), 
$$
and 
$$
\Delta_i(P_1,\ldots,P_g) = \det\big( (x^{(j)}_\al)_{\al = 0,\ldots, n,\al\neq i, 
j = 1,\ldots,n} \big) ,  
$$
where we set 
$$
x_\al^{(j)} = X_\al^{(j)} / X_0^{(j)} 
$$
(so $(x_\al^{(j)})_{\al = 1,\ldots,n}$ are the affine coordinates of 
$P_j$). 

Let $U \subset (\CC P^g)^{(g)}$ be the Zariski open subset defined 
as $\{(P_1,\ldots,P_g) \¦ X_0^{(1)}\neq 0,\ldots,X_0^{(g)} \neq 0, 
\Delta_0(P_1,\ldots,P_g) \neq 0\}$, then we get a map 
$$
(\CC P^g)^{(g)} \supset U \stackrel{(h_1,\ldots,h_g)}{\to} \CC^g. 
$$  
It turns out that $\varphi^{(g)}$ maps the generic point of $S^{(g)}$ to
$U$, so we get a map 
$$
S^{(g)} \supset V \stackrel{(h_1,\ldots,h_g) \circ \varphi^{(g)}}{\longrightarrow}
\CC^g. 
$$
This map is Beauville's fibration (\cite{Beauville}). 
In \cite{Mukai}, Mukai defined a Poisson structure on $S^{(g)}$,
which coincides with the symmetric power of the Poisson structure of $S$
on the smooth part of $S^{(g)}$. According to Proposition 3 of \cite{Beauville}, 
we have 

\begin{prop} \label{prop:Beauville}
$(h_1,\ldots,h_g) \circ \varphi^{(g)}$ is a lagrangian fibration. 
\end{prop}

This result may be derived from Theorem \ref{thm:Poisson} as follows: 
Proposition \ref{prop:Beauville} means that for any $i,j$, we have 
\begin{equation} \label{PB:comm}
\{(\varphi^{(g)})^*(h_i),  (\varphi^{(g)})^*(h_j)\} = 0 . 
\end{equation}
Since $\varphi^{(g)}(h_i) = \Delta_i / \Delta_0$, where 
$$
\Delta_i = [1,\varphi^*(h_1),\ldots,\check{\varphi^*(h_i)},\ldots,\varphi^*(h_n)], 
$$
and 
$$
\Delta_0 = [\varphi^*(h_1),\ldots,\varphi^*(h_n)], 
$$
(\ref{PB:comm}) follows from Theorem \ref{thm:Poisson}. 

\subsection{An affine version of Beauville's fibration} 
\label{affine:Beauville}

Let $S$ be a complex surface with Poisson structure, and let 
$\varphi : S \to \CC^g$ be an embedding. Set 
$$
\Delta_i = [1,\varphi^*(x_1),\ldots,\check{\varphi^*(x_i)},
\ldots,\varphi^*(x_g)], 
$$
for $i = 1,\ldots,g$, and 
$$
\Delta_0 = [\varphi^*(x_1),\ldots,\varphi^*(x_g)]. 
$$
Then if $\Delta_0$ is not zero, then the rational functions 
$$
h_i = \Delta_i / \Delta_0 
$$ 
on $S^{(g)} = S^g / \SG_g$, are Poisson commutative. 

\subsection{Beauville fibration in the case of the canonical cone of a 
curve} \label{2:4}

Let $C$ be an algebraic curve of genus $>1$. Let $K,\cO$ be its canonical
and structure sheaves. Set $\wt S_C = \PP(\cO\oplus K)$: $\wt S_C$ is 
a ruled surface, obtained from 
the total space of the cotangent bundle $T^* C$ by adding in each 
fibre, a point at infinity. We can blow down this additional copy of 
$C$ at infinity to a point. Let $S_C$ be  the resulting surface. 
The zero-section of $T^*C$ yields an embedding $C\subset S_C$. 
Let us define $\Gamma(S_C,\cO_{S_C}(* C))$ as the algebra of all 
rational functions on $S_C$, with only poles at $C$. 
Then we have an isomorphism of algebras 
\begin{equation} \label{isom}
\Gamma(S_C,\cO_{S_C}(* C)) = \bigoplus_{i\geq 0} H^0(C,K^{\otimes i}). 
\end{equation}
Indeed, an element of $H^0(C,K^{\otimes i})$ can be viewed as a rational 
function in each fiber of $T^* C \to C$, rational of degree
$-i$, and therefore as a function on each fiber of $\wt S_C \to C$, 
vanishing at $\infty$ and with only pole at $0$. 

The right side of (\ref{isom}) is the function algebra on the 
canonical cone $\Cone(C)$. Then (\ref{isom}) is also an isomorphism
of Poisson algebras: the Poisson structure of $\Gamma(S_C,\cO_{S_C}(* C))$ is 
induced by the symplectic structure of $T^* C$, and the Poisson structure of
$A_C = \bigoplus_{i\geq 0} H^0(C,K^{\otimes i})$ was defined in \cite{EO}. 

Then 
\begin{align*}
\¦nC\¦ & = \{f\in \Gamma(S_C,\cO_{S_C}(* C)) \¦ \on{val}_C(f)\geq -n\}
\\& = \oplus_{i=0}^n H^0(C,K^{\otimes i}). 
\end{align*}
Then the projective embedding corresponding to the linear system $nC$
is $S_C \hookrightarrow \PP(\¦nC\¦)$. 

On the other hand, we have an embedding 
\begin{equation} \label{embedding}
\Cone(C) \hookrightarrow \bigoplus_{i=1}^n H^0(C,K^{\otimes i}), 
\end{equation}
such that the diagram 
$$
\begin{array}{ccc}
S_C & \hookrightarrow & \PP(\¦nC\¦)
\\ 
\uparrow & \; & \uparrow
\\ 
\Cone(C) & \hookrightarrow  & \bigoplus_{i=1}^n H^0(C,K^{\otimes i})
\end{array}
$$
commutes. 

Let us construct the map (\ref{embedding}). 
Recall that the algebra of functions on $\Cone(C)$ is 
$$
A_C = \bigoplus_{i\geq 0} H^0(C,K^{\otimes i}). 
$$
The injection of 
vector spaces $\oplus_{i=1}^n H^0(C,K^{\otimes i}) \hookrightarrow 
A_C$ induces a morphism of algebras 
$$
S^\bullet(\oplus_{i=1}^n H^0(C,K^{\otimes i}) ) \to A_C, 
$$
dual to the map (\ref{embedding}). 

So the Beauville system associated to $(S_C,nC)$ corresponds to 
the system defined in Section \ref{affine:Beauville}, with respect to  
the embedding (\ref{embedding}). 

We now introduce a generalization of this system. 
Let $(d_1,\ldots,d_r)$ be integers $\geq 1$. Then the injection of 
vector spaces $\oplus_{i=1}^r H^0(C,K^{\otimes d_i}) \hookrightarrow 
A_C$ induces a morphism of algebras 
$$
S^\bullet(\oplus_{i=1}^r H^0(C,K^{\otimes d_i}) ) \to A_C, 
$$
and therefore an embedding 
$$
C \hookrightarrow \PP_{\on{weight}} \big(
\bigoplus_{i=1}^r H^0(C,K^{\otimes d_i})^* \big) ,  
$$
where $\PP_{\on{weight}}$ is the weighted projective space 
corresponding to the action of $\CC^\times$ given by 
$\la \cdot (v_1,\ldots,v_r) = (\la^{d_1}v_1,\ldots,\la^{d_r}v_r)$. 

Then the integrable system of Section \ref{affine:Beauville}
is defined by the Hamiltonians $h_i = \Delta_i / \Delta_0$, where 
we set $x_i = \omega_i$, and $(\omega_1,\ldots,\omega_N)$ is a basis of 
$\oplus_{i=1}^r H^0(C,K^{\otimes d_i})$. 

Let us give a direct proof of the Poisson commutativity of the 
$(h_i)_{i = 1,\ldots,N}$ in the particular case $r=1$. 

The Poisson structure on the function algebra $A_C
= \bigoplus_{i\geq 0} H^0(C,K^{\otimes i})$
may be defined as follows (see \cite{EO}): for any
rational form $\al$ on $C$, and any $i$-differential $\omega$, set 
$$
\nabla^\al(\omega) = \al^i d(\omega / \al^i),   
$$
and for any $i'$-differential $\omega'$, set 
$\{\omega,\omega'\} = i \omega \nabla^\al(\omega') - 
i' \omega' \nabla^\al(\omega)$. One checks that this definition is
independent of $\al$ and defines a Poisson structure on $A_C$. 

Then we get 
\begin{equation} \label{reason}
\{\Delta_i,\Delta_j\} = d_1 \cdot (\Delta_i F_j - \Delta_j F_i) , 
\end{equation}
where
$$
F_i = \sum_{j\neq i} [1, \omega_1,\ldots, \check \omega_i,\ldots, 
\nabla^\al(\omega_j), \ldots, \omega_N], 
$$
and 
$$
F_i = \sum_{j = 1}^N [\omega_1,\ldots, \nabla^\al(\omega_j), 
\ldots, \omega_g] . 
$$
Then (\ref{reason}) immediately implies $\{\Delta_i/\Delta_0,
\Delta_j/\Delta_0\} = 0$ for any $i,j$. 

\begin{remark} 
When $d_i = i$ for $i = 1,\ldots,r$, we obtain the "separated" 
version of Hitchin's system (see \cite{Hitchin,Donagi}).  
\end{remark}

\section{Quantization of Beauville fibrations}

The purpose of this section is to give a partial solution of the 
problem of quantizing the Beauville fibrations. In the 
next section, we will give a more explicit solution in the case of 
canonical cones.

\subsection{Quantization of fields}

Let $A$ be an integral algebra with Poisson structure, and let $K$  
be its fraction field. Let $A_\hbar$ be a quantization of $A$, i.e., 
$A_\hbar$ is a topologically free $\kk[[\hbar]]$-module, 
whose associated Poisson algebra is $A$. On the other hand, according to 
Lemma \ref{prol:poisson}, $K$ has a uniquely defined Poisson 
structure, extending the Poisson structure of $A$. 

\begin{prop} \label{quant:fields}
There is a unique quantization $K_\hbar$ of the Poisson 
ring $K$, containing $A_\hbar$ as a subalgebra. $K_\hbar$
is a skew field.  
\end{prop}

{\em Proof.} Let us first select a nonzero element $f_0$ of 
$A$ and construct the quantization of $A_{f_0}$. Let us fix 
$f\in A_\hbar$, whose reduction modulo $\hbar$ is $f_0$. 

We define $(A_\hbar)_f$ as the $\hbar$-adic completion of the quotient 
$(A_\hbar)[X] / I$, where $I$ is the vector space spanned by all 
$gf \otimes X^{n+1} - g \otimes X^n$, $g\in A_\hbar$, $n\geq 0$. 
The product is induced by the formulas 
$$
(aX^n) (b X^m) = \sum_{\al\geq 0} \pmatrix -n \\ \al
\endpmatrix a \ad(f)^\al(b) X^{n+m+\al}, 
$$
where the r.h.s. is $\hbar$-adically convergent. 
One checks that $(A_\hbar)_f$ is a quantization of $A_{f_0}$, 
and it is independent of the choice of $f$ above $f_0$. 
Repeating this construction for all nonzero elements of $A_0$, 
we construct $K_\hbar$. 
\hfill \qed \medskip

\subsection{Quantization of commuting families}

Assume that we are given a Poisson algebra $A$, 
such that $A^{\otimes n}$ is integral, and linearly
independent elements $f_0,\ldots,f_n$ $\in A$. To these data 
is associated a Poisson commuting family of elements $(h_i)_{i = 1,\ldots,n}$
of $\Frac(A^{\otimes n})$. By a quantization of the commuting 
family $(h_i)_{i = 1,\ldots,n}$, we understand: 

(1) a quantization $K_\hbar$ of the field $K = \Frac(A^{\otimes n})$

(2) a family of commuting elements  $(h_i)_\hbar$ of $K_\hbar$, 
deforming $h_i, i = 1,\ldots,n$. 

We will show: 

\begin{prop}
To construct a quantization of the commuting  family 
$(h_i)_{i = 1,\ldots,n}$, it suffices to construct a quantization of 
the Poisson algebra $A$.  
\end{prop}

Indeed, according to Proposition \ref{quant:fields}, the fraction field of 
$(A_\hbar)^{\otimes n}$ is a quantization of the fraction field of 
$A^{\otimes n}$; we then apply Theorem \ref{thm:comm}. 

\subsection{Quantization of Beauville systems}

To quantize the commuting families underlying the Beauville 
fibrations, it is therefore sufficient to quantize the coordinate rings of 
K3 surfaces. The solution of this problem is not known explicitly, in 
general. However, when $S$ is the canonical cone of an algebraic 
curve, a quantization is known in terms of formal pseudodifferential
operators (see \cite{EO}). Using the results of \cite{EO}, we can 
therefore quantize the Beauville systems in this case. In the next section, 
we will make this solution explicit.

\section{The case of the canonical cone of an algebraic curve}
\label{sect:final}

Recall the situation of Section \ref{2:4}. The surface $S$ is 
birationally equivalent to the canonical cone $\Cone(C)$ of 
an algebraic curve $C$, and we have an embedding
$$
\Cone(C) \hookrightarrow \bigoplus_{i=1}^r H^0(C,K^{\otimes d_i})^*. 
$$
To these embeddings correspond classical integrable systems. 
We explained how to construct their quantizations. We will 
show that when $r=1$, these quantized integrable systems can be 
obtained as a commuting family of differential operators in 
symmetric powers of $C$. 

\subsection{Algebras of rational differential operators}

Let $N$ be an integer, and let $\DO_\rat(C^N)$ be the 
algebra of rational differential operators on $C^N$. 
If $X$ is a fixed nonzreo rational vector field on $C$, 
then $\DO_\rat(C^N)$ is a subalgebra of $\End(\CC(C^N))$, 
and an element of this algebra is uniquely written as 
\begin{equation} \label{op:type}
\sum_{\al_1,\ldots,\al_N} f_{\al_1,\ldots,\al_N}
X^{\al_1} \otimes \cdots \otimes X^{\al_N},   
\end{equation}
where $f_{\al_1,\ldots,\al_N} \in \CC(C^N)$ and all 
but finitely many $f_{\al_1,\ldots,\al_N}$ are zero. 

Set $\Fil_i(\DO_\rat(C^N)) = \{$operators of the form (\ref{op:type}), 
such that $f_{\al_1,\ldots,\al_N} = 0$ when $\al_1 + \cdots + \al_N >i\}$. 
This defines an algebra filtration on $\DO_\rat(C^N)$. The associated graded 
algebra identifies the  algebra 
$$
\CC(C^N)[\xi_1,\ldots,\xi_N]
$$
(the tensor product of $\CC(C^N)$ with a polynomial algebra). 

\subsection{Relation with $\Cone(C)$}

Let $i$ be any integer. The space of all rational functions on $\Cone(C)$, 
homogeneous of degree $i$ along the fibers of $\Cone(C)\to C$, 
identifies with 
$$
\{\on{rational\ }i\on{-differentials\  on\ }C\} = 
\{\on{rational\ sections\ of\ }K^{\otimes i}\}. 
$$
The direct sum 
$$
\bigoplus_{i\in\ZZ}
\{\on{rational\ sections\ of\ }K^{\otimes i}\}
$$
is a subalgebra of $\CC(\Cone(C)) = \CC(S_C)$. Moreover, 
there is a unique algebra morphism 
$$
\CC(C)[\xi] \to 
\bigoplus_{i\in\ZZ}
\{\on{rational\ sections\ of\ }K^{\otimes i}\}
\subset \CC(\Cone(C)), 
$$
taking each $f\xi^n$ to $fX^n$ (a rational section of $K^{\otimes -n}$). 

In the same way, for $(i_1,\ldots,i_N)$ a sequence of integers, 
the space of rational functions on $\Cone(C^N)$, homogeneous of 
degree $(i_1,\ldots,i_N)$ in the fibres of $\Cone(C)^N \to C^N$, 
is 
$$
\{\on{rational\ sections\ of\ }K^{\otimes i_1} \boxtimes \cdots \boxtimes 
K^{\otimes i_N} \on{\ over\ }C^N\}, 
$$
and we have an algebra morphism 
$$
\CC(C^N)[\xi_1,\ldots,\xi_N] \to \bigoplus_{(i_1,\ldots,i_N)\in\ZZ^N}
\{\on{rational\ sections\ of\ }K^{\otimes i_1} \boxtimes \cdots \boxtimes 
K^{\otimes i_N}\} \subset \CC(\Cone(C)^N).  
$$

\subsection{Commuting differential operators}

Let us assume that we are in the situation of Section \ref{2:4}, 
and that $r=1$. We set $N = \dimm(H^0(C,K^{\otimes d_1}))$. 

\begin{lemma} When $i = 1,\ldots,N$, we have 
$$
H_i^\cl \in \CC(C^N)[\xi_1,\ldots,\xi_N]. 
$$
\end{lemma}

{\em Proof.} We have 
$$
\Delta_i^\cl = \sum_{j=1}^N f_{i,j} 
(\xi_1\cdots \check \xi_i \cdots \xi_N)^{-d}, 
$$
and 
$$
\Delta_0^\cl = f
(\xi_1\cdots \xi_N)^{-d}, 
$$
where $f,f_{i,j}$ belong to $\CC(C^N)$. So 
$H_i^\cl = \sum_{j=1}^N (f_{i,j}/f) (\xi_j)^d$. 
\hfill \qed \medskip 

\begin{prop}
There exists a commuting family $(H_1,\ldots,H_N)$
of rational differential operators on $C^N$, with symbols 
$(H_1^\cl,\ldots,H_N^\cl)$. 
\end{prop}

{\em Proof.}
Let $(\omega_1,\ldots,\omega_N)$ be a basis of 
$H^0(C,K^{\otimes d_1})$ and let $X$ be a nonzero rational 
vector field on $C$. Set $f_i = \omega_i X^d$ (product of sections of 
bundles), then each $f_i$ is a rational function on $C$. Moreover, $\wt\omega_i =
f_i X^{-d}$ is a formal pseudodifferential operator on $C$, with symbol 
$\omega_i$ (see \cite{EO}). Let us compute the $\Delta_i$ and $H_i$
corresponding to the family $(\wt\omega_1,\ldots,\wt\omega_N)$. 
We have 
$$
\Delta_0 = \Phi (X^{-d} \otimes \cdots \otimes X^{-d}), 
$$
where 
$$
\Phi = \sum_{\sigma\in\SG_N} \eps(\sigma) 
f_{\sigma(1)} \otimes \cdots \otimes f_{\sigma(N)} ,  
$$
and
$$
\Delta_i = \sum_{i=1}^N (-1)^{Ni} \delta_i^{(1\ldots \check i \ldots N)}, 
$$
where 
$$
\delta_i^{(1\ldots N-1)} = \Phi_i  (X^{-d} \otimes \cdots \otimes X^{-d}), 
$$
and 
$$
\Phi_i = \sum_{\sigma\in\SG_{N-1}} \eps(\sigma)
g^{(i)}_{\sigma(1)} \otimes \cdots \otimes g^{(i)}_{\sigma(N-1)},  
$$
and $(g^{(i)}_1,\ldots,g^{(i)}_N) = (f_1,\ldots,\check f_i, \ldots, f_N)$. 
So we get 
$$
H_i = \Delta_i (\Delta_0)^{-1} 
= \sum_{i=1}^N h_i^{(i1\ldots \check i \ldots N)}, 
$$
where
$$
h_i^{(1\ldots N)} = \Phi_i (X^d \otimes 1 \otimes \cdots \otimes 1) 
\Phi^{-1} .  
$$
So $h_i^{(1\ldots N)}$ is a rational differential operator, and so is $H_i$. 
\hfill \qed \medskip 

\subsection{Explicit formulas in the rational case}

Applying Theorem \ref{thm:comm} to the family 
$$
f_i = T^{-1} \cdot {1\over{z-P_i}}, \; i = 1,\ldots,N, 
$$
for $T$ any rational differential operator on $\CC P^1$, we get: 

\begin{thm}
Define $T_{z_i}$ as the differential operator on 
$(\CC P^1)^N$, acting as $T$ on the $i$th variable. 
Let $(P_1,\ldots,P_N)$ be a set of distinct points of $\CC$. Set 
$$
H_k = \sum_{i=1}^N 
{{\prod_{(i',k')\¦ i' = i \on{\ or\ }k'=k} (z_{i'} - P_{k'}) }
\over{\prod_{i'\¦ i'\neq i} (z_i - z_{i'})}} \cdot T_{z_i}
$$ 
for $k = 1,\ldots,N$. Then $(H_1,\ldots,H_N)$ is a 
commuting family of rational differential operators on 
$(\CC P^1)^N$. 
\end{thm}

Indeed, we have $H_k = \Delta_k (\Delta_0)^{-1}$.  

\appendix
\section{Relation of Theorem \ref{thm:comm} with formal noncommutative geometry}
\label{app}

In this section, we assume $\on{char}(\kk) = 0$. 

\subsection{}

Let $x_1,\ldots,x_n$ be formal variables, and let $\on{FreeAlg}(x_1,\ldots,x_n)$
be the free algebra with generators $x_1,\ldots,x_n$. This is the enveloping algebra of
the free Lie algebra with the same generators, $\on{FreeLie}(x_1,\ldots,x_n)$, so 
symmetrization induces a linear isomorphism 
\begin{equation} \label{isom:free}
S^\bullet(\on{FreeLie}(x_1,\ldots,x_n)) \to 
\on{FreeAlg}(x_1,\ldots,x_n). 
\end{equation}  
Moreover, we can define a grading on $S^\bullet(\on{FreeLie}(x_1,\ldots,x_n))$
by giving degree $k-1$ to an element of $\on{FreeLie}(x_1,\ldots,x_n)$ of 
degree $k$. Then the algebra structure of $S^\bullet(\on{FreeLie}(x_1,\ldots,x_n))$
induced by (\ref{isom:free}) extends uniquely to its completion for this 
grading. We denote by $\wh{\on{FreeAlg}}(x_1,\ldots,x_n)$ the resulting 
completed algebra. 

Let us denote by $\on{FreePoisson}(x_1,\ldots,x_n)$ the free Poisson 
algebra with generators $x_1,\ldots,x_n$. Then 
$\on{FreePoisson}(x_1,\ldots,x_n)$ is isomorphic to 
the symmetric algebra $S^\bullet(\on{FreeLie}(x_1,\ldots,x_n))$. 
Then we have 
$$
\on{gr}(\on{FreeAlg}(x_1,\ldots,x_n)) = 
\on{FreePoisson}(x_1,\ldots,x_n) .  
$$
We denote by $\wh{\on{FreePoisson}}(x_1,\ldots,x_n)$ the completion of 
$\on{FreePoisson}(x_1,\ldots,x_n)$ for the same grading as above. 

\subsection{} 

Let $F_n$ be the algebra with generators $f_{i,k}$, 
$i = 0,\ldots,n$, $k = 1,\ldots,n$, and relations 
$f_{i,k} f_{j,\ell} = f_{j,\ell} f_{i,k}$ when $k\neq \ell$. 
Then $F_n$ is isomorphic to the tensor product 
$\otimes_{k=1}^n \on{FreeAlg}(f_{0,k},\ldots,f_{n,k})$, 
where the tensor factors commute with each other. 
Set 
$$
\wh F_n = \bigotimes_{k=1}^n \wh{\on{FreeAlg}}(f_{0,k},\ldots,f_{n,k}) . 
$$
Set 
$$
\Delta_0 = \sum_{\sigma\in\SG_n} \eps(\sigma) 
f_{1,\sigma(1)} \cdots f_{n,\sigma(n)} . 
$$
According to \cite{Kap}, we can localize $\wh F_n$ with respect
to $\Delta_0$. 

Moreover, define $P_n$ as the Poisson algebra 
$\otimes_{k=1}^n \on{FreePoisson}(f_{0,k},\ldots,f_{n,k})$, 
where the tensor factor Poisson commute with each other, and 
$\wh P_n$ as its completion 
$$
\wh P_n = \bigotimes_{k=1}^n \wh{\on{FreePoisson}}(f_{0,k},\ldots,f_{n,k}) . 
$$
Let $\Delta_0^{\on{Poisson}}$ be the analogue of $\Delta_0$
in $\wh P_n$, then we can localize $\wh P_n$ with respect to
$\Delta_0^{\on{Poisson}}$. Moreover, we have 
$$
\gr((\wh F_n)_{\Delta_0}) = (\wh P_n)_{\Delta^{\on{Poisson}}_0}. 
$$ 

Then 
\begin{prop}
Set $\Delta_i =  \sum_{\sigma\in\SG_n} \eps(\sigma) 
f_{0,\sigma(1)} \cdots f_{i-1,\sigma(i)} f_{i+1,\sigma(i+1)} 
\cdots f_{n,\sigma(n)}$. Then the elements $H_i = \Delta_i (\Delta_0)^{-1}$
commute with each other. 
\end{prop}

{\em Proof.} Let us denote by $(\delta_\al)_\al$ the collection of 
all minors obtained from the family $(f_{i,k})_{0\leq i\leq n,1\leq k\leq n}$. 
Then we can also localize $\wh F_n$ with respect to this family, we we have 
a sequence of inclusions 
$$
\wh F_n \hookrightarrow (\wh F_n)_{\Delta_0}
\hookrightarrow (\wh F_n)_{\Delta_0,(\delta_\al)}. 
$$
According to Corollary \ref{comm:fij}, the images of the $H_i$ in the last algebra
commute together. This implies that the $H_i$ already commute in 
$(\wh F_n)_{\Delta_0}$. 
\hfill \qed \medskip 

We also get the following result. 
\begin{prop}
Let $A$ be an algebra. Let us define $\on{Comm}_n(A)$ as the linear span 
of all order $n$ commutators, and  $\on{Fil}^i(A)$ as the sum 
$$
\sum_{n_1,\ldots,n_k\¦ n_1+ \cdots + n_k \geq i} \on{Comm}_{i_1}(A)
\cdots \on{Comm}_{i_k}(A).
$$
Then in Corollary \ref{comm:fij}, the hypothesis ``the $\Delta_{I,J}$ are 
all invertible" may be replaced by ``$\Delta_0$ is invertible
and $A$ is complete and separated for the topology defined by $\on{Fil}^i(A)$". 
\end{prop}

Indeed, the hypothesis implies that we have an algebra morphism 
$(\wh F_n)_{\Delta_0} \to A$. 

\subsection{}

Let us now discuss the deformation of the commuting 
family $(H_i)_{i=1,\ldots,n}$.

Let us set 
$$
T_{\on{Poisson}} = \{(h_1,\ldots,h_n) \¦ h_1,\ldots,h_n\in 
(\wh P_n)_{\Delta_0^{\on{Poisson}}} \; \on{and}\; 
\{H^{\on{Poisson}}_i,h_j\}+\{h_i,H^{\on{Poisson}}_j\} = 0 \}
$$
and 
$$
T_{\on{assoc}} = \{(h_1,\ldots,h_n) \¦ h_1,\ldots,h_n\in 
(\wh F_n)_{\Delta_0} \; \on{and}\; 
[H_i,h_j]+[h_i,H_j] = 0 \}. 
$$
Then $T_{\on{Poisson}}$ contains the families 
$h_i = \sum_j \la_{i,j} H^{\on{Poisson}}_j$ and $h_j = \{a,H^{\on{Poisson}}_i\}$, 
for $\la_{i,j}\in\kk$ and $a\in (\wh P_n)_{\Delta_0^{\on{Poisson}}}$. 
In the same way, $T_{\on{assoc}}$ contains the families 
$h_i = \sum_j \la_{i,j} H_j$ and $h_j = [a,H_i]$, 
for $\la_{i,j}\in\kk$ and $a\in (\wh F_n)_{\Delta_0}$. 

Equality of  $T_{\on{Poisson}}$ (resp., $T_{\on{assoc}}$)  
with its subspace means that the commuting family $(H_i)_{i=1,\ldots,n}$ 
(resp., $(H^{\on{Poisson}}_i)_{i=1,\ldots,n}$) has no
nontrivial deformations. It is easy to check that the absence of 
nontrivial deformations in the Poisson situation implies the same 
statement in the associative situation. 

\subsection*{Acknowledgements} 

{\small
We would like to thank A. Odesskii and M. Olshanetsky for discussions
on the subject of this work. We would also like to thank the Mathematisches Institut
Oberwolfach for hospitality at the time this work was done. V.R. was 
partially supported by grants INTAS 99-1705, RFBR 01-01-00549 and by the grant 
for scientific schools RFBR 00-15-96557. }

\end{document}